\documentclass[a4j,12pt]{article}
\usepackage{amsmath,amsthm,amssymb,bm}

\numberwithin{equation}{section}

\newtheorem{thm}{Theorem}[section]

\newtheorem{cor}[thm]{Corollary}
\newtheorem{definition}[thm]{Definition}
\newtheorem{remark}[thm]{Remark}
\newtheorem{example}[thm]{Example}

\newenvironment{df}{\begin{definition}\rm}{\end{definition}}
\newenvironment{rem}{\begin{remark}\rm}{\end{remark}}
\newenvironment{ex}{\begin{example}\rm}{\end{example}}

\title{Operational extreme points of unital completely positive maps}

\author{Marie Choda}

\date{}

\begin{document}
\maketitle
\centerline{Osaka Kyoiku University, Osaka, Japan} 
\centerline{marie@cc.osaka-kyoiku.ac.jp}
\begin{abstract} 
Two notions for linear maps (operational convex combinations and operational 
exetreme points)   are introduced. 
The set $S$ of ucp maps on $M_n(\mathbb {C})$ is the 
operational convex combinations of the identity map. 
An operational extreme point of $S$ is an extreme point of $S$ but the converse does not hold, and 
every automorphism is an operational extreme point of $S$. 
\end{abstract} 

keywords: {Positive linear map, convex combination}

\smallskip

\section{Introduction} 
From a view point of von Neumann entropy for states of $M_n(\mathbb {C})$, 
we gave some characterization for unital positive   Tr-preserving  maps of the 
algebra of $n \times n$ conplex matrices $M_n(\mathbb {C})$. 
That is, a positive unital  Tr-preserving map $\Phi$ of $M_n(\mathbb {C})$ preserves the von Neumann entropy 
of a given state $\phi$ if and only if $\Phi$ plays a role  of an automorphism for $\phi$.  

In this note, we pick up the set of  unital completely positive (called "ucp" for short) maps of $M_n(\mathbb {C})$. 
That is, for a unital linear map of $M_n(\mathbb {C})$, 
we replace the property "Tr-preserving and positive "   to the property "completely positive", 
and investigate that what kind of position the automorphisms stand in ucp maps. 
\vskip 0.3cm

The set of ucp maps is a convex subset of linear maps of $M_n(\mathbb {C})$. 
The notion of convex sets begins with basic definition of the linear concepts of addition and scalar 
multiplication. 

Here, we shall consider the notion of convexity 
not only scalar multiplication but also the multiplication via operators and 
generalize the notion of convexity, i.e.,
we introduce the notion of operational convex conbination. 

The motivation for the terminology "operational convex conbinations" comes from 
the following two definitions: One is 
Lindblad's "operational partition"  in \cite{L} (cf. \cite{NS} or \cite{OP}) and the other is 
Cuntz's canonical endomorphism $\Phi_n$ in \cite{Cu}. 
It seems to be natural for treating the set of ucp maps as the set of operational convex conbinations of 
automorphisms of $M_n(\mathbb {C})$. 

We also introduce the notion of operational extreme point. 
Since an automorphism of $M_n(\mathbb{C})$ is an extreme point of the set of positive maps (cf. \cite{St}), 
any automorphism can not be expressed as a convex combination of  two different 
automorphisms.  
However, it is possible to be expressed as an operational convex combination of  two different automorphisms. 
By our definition, operational extreme points are extreme points, but 
the converse is not true as we show in Example 3.6, and 
we show  that automorphisms are  operational extreme points in the set of ucp  maps of  $M_n(\mathbb {C})$. 

\section{Preliminaries}
Here we summarize notations, terminologies and basic facts. 

\subsection{Finite partition of unity} 
What we need to define a convex sum ? 
In usual, we need a probability vector $\lambda = (\lambda_1, \cdots, \lambda_n )$:  
$\lambda_i \geq 0, \quad \sum_i  \lambda_i = 1.$ 

Given a finite subset $x = \{x_1, ..., x_n \}$ of  a vector space $X$, 
the vector  
$\sum_i  \lambda_i x_i$ is  called a convex sum of $x$ via $\lambda$. 

Now, we consider such a $\lambda$  as a "finite partition of 1". 
\vskip 0.3cm

Two  generalized notions of  finite partition of 1 
are given in  the framework of the non-commutative entropy as follows: 

Let $A$ be a unital $C^*$-algebra.  

(1) A finite subset $\{x_1, ..., x_k \}$ of  $A$ is called  a 
{\it finite partition of unity}  by Connes-St\o rmer (\cite {CS}) if 
they are nonnegative operators which satisfy that 
$\sum_{i=1}^n  x_i = 1_A$, 

(2) A finite subset $\{x_1, ..., x_k \}$ of  $A$ is  called  a 
{\it finite  operational partition in $A$ of unity of size $k$} 
by Lindblad (\cite{L})     if 
$\sum_i^k x_i^* x_i = 1_A$.  
\vskip 0.3cm

Our main target in this note is a finite subset $ \{v_1, ..., v_k \}$ of non-zero elements in $A$ such that 
$ \{v_1^*, ..., v_k^* \}$ is a finite  operational partition of unity so that   
$\sum_i^k v_i v_i^* = 1_A.$   
We call such a set $ \{v_1, ..., v_k \}$ a {\it finite operational partition of unity of size $k$} in $A$,  
and denote  the set of all  finite operational partition of unity in $A$ by  $FOP(A)$: 
\begin{equation}
FOP(A) = \{  \{v_1, ..., v_k \} \ | \ 0 \ne v_i \in A, \forall i, \ \sum_i^k v_i v_i^* = 1_A, \  k = 1, 2, \cdots \}
\end{equation}
We denote by $U(A)$ the set of all unitaries in $A$.  Clearly, $U(A)$  is the set of the most trivial finite operational partition of 
unity with the  size $1$. 

\subsection{Unital completely positive  (ucp) map  $\Phi$} 
Let $M_n(\mathbb {C})$ be the $C^*$-algebra  of $n \times n$ matrices over the complex field $\mathbb {C}$.
A linear map $\Phi$ on a unital $C^*$-algebra $A$ is  positive iff $\Phi(a)$ is  positive for all   positive $a \in A$ 
and {\it completely positive} iff $\Phi \otimes 1_k$ is  positive for all   positive integer $k$, where 
the map $\Phi \otimes 1_k$ is the map on $A \otimes M_k(\mathbb {C})$ defined by 
$\Phi \otimes 1_k (x \otimes y) = \Phi(x) \otimes  y$ for all $x \in A$ and $y \in M_k(\mathbb {C})$.  

We restrict the unital $C^*$-algebra $A$ to  $M_n(\mathbb {C})$. 

In \cite [Theorem 2] {Choi2}, Choi gave the following characterization:    
a linear map $\Phi$ of  $M_n(\mathbb {C})$ is completely positive iff
$\Phi$ is of the form $\Phi(x) = \sum_{i=1}^m v_i x v_i^*$ for all $x \in M_n(\mathbb {C})$ 
by some $\{v_i\}_{i=1}^m \subset M_n(\mathbb {C})$. 
Moreover,  
for $\{v_i\}_{i=1}^m$ inducing the form $\Phi(x) = \sum_{i=1}^m v_i x v_i^*$, 
we may require that $\{v_i\}_i$  is linearly independent 
so that in the form  the number $m$  is  uniquely determind. 
Such a form was called a 'canonical' expression for $\Phi$ (see \cite [Remark 4] {Choi2}). 

Let us call the  uniquely determind number $m$ the {\it size of the $\Phi$}.  
\vskip 0.3cm 

Now we pick up the case where $\Phi$ is a 
unital completely positive (called "ucp" for short)  map  of $M_n(\mathbb {C})$. 
Then the $\{v_1, \cdots, v_m\} \subset M_n(\mathbb {C})$ used in the form 
$\Phi(x) = \sum_{i=1}^m v_i x v_i^*,  (x \in M_n(\mathbb {C}))$ 
satisfies that $\sum_{i=1}^m v_i v_i^* = 1$. 
This means that each ucp map $\Phi$ of $M_n(\mathbb {C})$ is induced some $\{v_1, \cdots, v_m\}$ in $FOP(M_n(\mathbb {C}))$. 

Given an operator $v \in M_n(\mathbb {C})$, the map ${\rm Ad}v$ on $M_n(\mathbb {C})$ is given  by 
${\rm Ad}v (x) = vxv^*, (x \in M)$.  
Then the  group  $Aut(M_n(\mathbb{C}))$ of all automorphisms of $M_n(\mathbb{C})$ is 
written by the  form $Aut(M_n(\mathbb{C})) =  \{ {\rm Ad}u \ | \ u \in U(M_n(\mathbb{C})) \}$, 
where $U(M_n(\mathbb{C}))$ is the group of all unitaries in $M_n(\mathbb{C})$. 
Similarly,  the set $UCP(M_n(\mathbb{C}))$  of all ucp maps on $M_n(\mathbb{C})$ 
is written by the following form: 
\begin{equation}
UCP(M_n(\mathbb{C}))  = \{\sum_{i=1}^m {\rm Ad}v_i \ | \ \{v_i\}_{i=1}^m \in FOP(M_n(\mathbb{C})),  \ m=1,2,\cdots   \}
\end{equation}

\section{\bf Operational Convex Combination} 

\subsection{\bf Operational convexity} 

\begin{df} 
Let $\{\Phi_i \}_{i=1}^m$ be 
a set of linear maps  on $M_n(\mathbb{C})$ and $\{v_i\}_{i=1}^m \in FOP(M_n(\mathbb{C}))$. 
We call $\sum_{i=1}^m {\rm Ad}v_i \circ \Phi_i$ 
an {\it operational convex combination} of $\{\Phi_i \}_{i=1}^m$ with an {\it operational coefficients} $\{v_i\}_{i=1}^m$. 
We also say that a subset $S$  of linear maps  on $M_n(\mathbb{C})$ 
is  {\it operational convex} if it is closed under all operational convex combinations. 
\end{df}

We can consider 
$UCP(M_n(\mathbb{C}))$ as the set of all operator convex combinations of the group $Aut(M_n(\mathbb{C}))$. 
Moreover $UCP(M_n(\mathbb{C}))$ is represented as the set of all operational convex combinations of the identity $id$ of $M_n(\mathbb{C})$. 
We give some characterization for  a role of $Aut(M_n(\mathbb{C}))$ in $UCP(M_n(\mathbb{C}))$ from a view point of extreme points. 

\subsubsection{Cuntz's canonical endomorphism as an example} 
The Cuntz's canonical endomorphism $\Phi_n$ (\cite{Cu}) is an interesting example in unital completely positive maps of infinite 
dimensional simple $C^*$-algebras, which is given as an operational convex combination of the identity. 
That is, let $\{S_1, S_2, \cdots, S_n \}$ be isometries on an infinite dimensional Hilbert space $H$ such that 
$\sum_i S_i S_i^* = 1$. The Cuntz algebra $O_n$ is the $C^*$-algebra generated by $\{S_1, S_2, \cdots, S_n \}$. 
The map $\Phi_n$ is given as $\Phi_n (x) = \sum_i S_i x S_i^*$ for all $x \in O_n$. So, in our notation, 
$\{S_1, S_2, \cdots, S_n \} \in FOP(O_n)$  and $\Phi_n \in UCP(O_n)$.   

The left inverse $\Psi$ of $\Phi_n$ plays an inportant role in the theory of Cuntz algebras and it 
is  given by the form $\Psi(x) = (1/n) \sum_i S_i^* x S_i, (x \in O_n)$. 

We remark that $\Psi$ is also an  operational convex combination of the identity and 
$\Psi \in UCP(O_n)$. 

Later we discuss in another paper on the case of unital infinite dimensional $C^*$-algebras represented by $O_n$. 

\subsubsection{Operational extreme point} 
Now let us remember the  notion of  extreme points. 
Let $S$ be a convex set. 
Then a $z \in S$ is an  extreme point in $S$ if $z$ cannot be the convex combination 
$\lambda x + (1 - \lambda) y$ of two points $x, y \in S$ with $x \ne y$ and $\lambda \in (0,1)$, 
i.e., if $z = \lambda x + (1 - \lambda) y, (x, y \in S)$ then $x = y = z$. 

In this note, we say this notion of extreme points   an extreme point {\it in the usual sense}. 

\begin{rem}
A related notion are investigated for positive maps on $C^*$-algebras in \cite{St}. 
A positive map $\Phi$ on a $C^*$-algebra $A$  is  {\it extremal} if the only positive maps 
$\Psi$ on $A$, such that $\Phi - \Psi$ is positive, are of the form $\lambda \Phi$ with $0 \le \lambda \le 1$.  
In the set of all positive maps on $B(H)$ for a Hilbert space $H$, the map  ${\rm Ad} u, (u \in B(H))$ is extremal 
\cite [Proposition 3.1.3] {St}. 
This implies that 
any automorphism of $M_n(\mathbb{C})$ can not be expressed as a convex combination of  two different 
automorphisms. 
However if we replace a convex combination to an oparational convex combination, then it is possible 
for each automorphism $\Phi$ of $M_n(\mathbb{C})$.    
\end{rem}

\begin{ex} 
Let $\Theta, \Phi$ and $\Psi$ be three different automorphisms of $M_n(\mathbb{C})$. 
Assume that $u, v$ and $w$ are unitaries in $M_n(\mathbb{C})$ such that $\Theta = {\rm Ad} u, 
\Phi = {\rm Ad} v$ and $\Psi = {\rm Ad} w$. 
Put $a = \lambda^{1/2} uv^*$ and  $b = (1-\lambda)^{1/2} uw^*$ for some $\lambda \in (0,1)$. 
Then  $\{a, b\} \in FOP(M_n(\mathbb{C}))$, and 
the operational convex conbination of  $\Phi$ and $\Psi$ 
with the operational coefficients $\{a, b\}$ is the automorphism $\Theta$, i.e., 
$a \Phi(x){a}^* +  b \Psi(x) {b}^* = \Theta(x)$ for all $x \in M_n(\mathbb{C})$. 
\end{ex}
\vskip 0.3cm

In the case of operational convex combinations for linear maps $\Phi$ and $\Psi$ on $M_n(\mathbb{C})$ with 
an operational coefficient $\{a, b\} \in FOP(M_n(\mathbb{C})$, the map 
${\rm Ad}a \circ \Phi$ corresponds $\lambda x$ and the ${\rm Ad}b \circ \Psi$ does $(1 - \lambda) y$. 
Putting this in mind, let us  define as  follows and show that an automorphism of $M_n(\mathbb{C})$ (i.e., the ucp maps with the size $1$) 
is an operational extreme point. 
\vskip 0.3cm

\begin{df} 
Let $S$ be an operational convex subset of linear maps on $M_n(\mathbb{C})$. 
We say that a $\Phi \in S$ is an {\it operational extreme point} of $S$  
if a reprensentation of $\Phi$ that 
$\Phi = {\rm Ad}a \circ \Phi_1 + {\rm Ad}b \circ \Phi_2, (\Phi_i \in S,  (i = 1,2), \{a, b\} \in FOP(M_n(\mathbb{C}))$ 
implies  that 
$aa^* = \lambda 1_{M_n(\mathbb{C})}, \ bb^* = (1 - \lambda) 1_{M_n(\mathbb{C})}$ for some $\lambda \in (0,1)$ 
so that  
$\lambda^{-1} {\rm Ad}a \circ \Phi_1  = \{1-\lambda \}^{-1} {\rm Ad}b \circ \Phi_2 = \Phi$. 
\end{df}

\begin{rem}
If $\Phi$ is an  operational extreme point of a  convex subset $S$  of linear maps on $M_n(\mathbb{C})$, 
then  $\Phi$ is an  extreme point of $S$ in the usual sense.  
\end{rem}

In fact if $\Phi = \lambda \Phi_1 + (1 - \lambda) \Phi_2, (\Phi_i \in S)$ and if 
$\Phi$ is an  operational extreme point of $S$. 
then $\Phi_1  = \lambda^{-1}  \lambda \Phi_1  = \Phi$ and $\Phi_2 = \{1-\lambda \}^{-1} (1-\lambda ) \Phi_2 = \Phi$ 
so that $\Phi$ is an  extreme point of $S$ in the usual sense. 
\vskip 0.3cm

The converse of the above is not true in general as the following example shows: 

\begin{ex} Here we give a ucp map which is an extreme but not operational extreme point.   

Let $\{e_{ij}\}_{i, j = 1,2}$ be a matrix units of $M = M_n(\mathbb{C})$, and let 
$V = \{v_1 = e_{11}, v_2 = e_{21}\} \in FOP(M)$. 
Let $\Phi_V = \sum_{i = 1}^2 {\rm Ad} v_i$. 
Then the ucp map $\Phi_V$ satisfies that 
$\Phi_V(e_{11}) = 1_M$, and $\Phi_V(e_{12}) = \Phi_V(e_{21}) = \Phi_V(e_{22}) = 0$. 
Assume that $\Phi_V$ has a form  $\Phi_V = \lambda \Phi_1 + (1 - \lambda) \Phi_2$ by $\Phi_i \in UCP(M), (i = 1,2)$. 
Then $0 = \Phi_V(e_{22}) = \lambda \Phi_1(e_{22}) + (1 - \lambda) \Phi_2(e_{22})$. 
This implies that $\Phi_i(e_{22}) = 0, (i = 1,2)$ because $\Phi_i(e_{22})$ is positive for $i = 1,2$ 
and that  $\Phi_i(e_{11}) = \Phi_i(1_M - e_{22}) = 1_M  - 0 = 1_M $ for $i = 1,2$. 
By using Kadison-Schwarz inequality (cf. \cite{NS, OP}), we have that  
$\Phi_i(e_{21}) \Phi_i(e_{21})^* \leq  \Phi_i(e_{21} e_{21}^*) = \Phi_i(e_{22}) = 0$  so that 
$\Phi_i(e_{21}) = 0$ for $i = 1,2$. 

Hence $\Phi_V = \Phi_1 = \Phi_2$ and $\Phi_V$ is an  extreme point of $UCP(M)$.  
\smallskip

Now as a $\{a, b\} \in FOP(M)$ we choose  $a = e_{11}, b = e_{22}$. 
As two ucp maps $\Phi_3, \Phi_4$ on $M$, let $\Phi_3$ be the identity map  and $\Phi_4$ be 
the ${\rm Ad} u$ where $u = e_{12} + e_{21}.$  
Then $\Phi_V = {\rm Ad} a \circ \Phi_3 + {\rm Ad} b \circ \Phi_4$. 
This  shows that $\Phi_V$ is not an operational extreme point. 
\end{ex}

\begin{thm} 
If an automorphism $\Theta$ of $M_n(\mathbb{C})$ is decomposed into 
an operational convex combination of the form that 
$\Theta ={\rm Ad}a \circ \Phi  +{\rm Ad} b \circ \Psi$ 
via $\{a, b\} \in FOP(M_n(\mathbb{C}))$ and $\Phi, \Psi  \in UCP(M_n(\mathbb{C}))$, 
then there exist  unitaries $u_a, u_b \in M_n(\mathbb{C})$ and a $\lambda \in (0,1)$ 
such that 
\begin{equation}
a = \sqrt{\lambda} u_a, \ b = \sqrt{1 - \lambda} u_b \quad \text{so that} \quad 
\Phi = {\rm Ad} u_a^*u, \ \Psi = {\rm Ad} u_b^*u.
\end{equation} 
Here $u$ is a unitary  with $\Theta = {\rm Ad} u$. 
\end{thm}

{\it Proof}. 
Since $\Theta - {\rm Ad}a \circ \Phi (= {\rm Ad} b \circ \Psi)$ is positive and $\Theta$ is extremal ([Proposition 3.1.3, \cite{St}]), 
there exists a $\lambda \in [0, 1]$ such that ${\rm Ad}a \circ \Phi = \lambda {\rm Ad} u$, which implies that 
${\rm Ad}b \circ \Psi = ( 1 - \lambda) {\rm Ad} u$. 
On the other hand, since $\Phi$ and $\Psi$ are unital, we have that 
$aa^* = {\rm Ad}a \circ \Phi(1) = \lambda {\rm Ad} u (1) = \lambda 1$ and $bb^* = (1- \lambda) 1$. 
Hence we have unitaries  $u_a, u_b \in M_n(\mathbb{C})$ such that $a = \sqrt{\lambda} u_a$ and  $b = \sqrt{1 - \lambda} u_b$.  
These relations imply that 
$\lambda u_a \Phi(x) u_a^* = {\rm Ad}a \circ \Phi(x) = \lambda u x u^*$ so that $\Phi = {\rm Ad} u_a^*u$. 
Similarly, $\Psi = {\rm Ad} u_b^*u$. 
\qed
\vskip 0.3cm

\begin{cor} 
An automorphism on $M_n(\mathbb {C})$ is an operational extreme point of the set of 
the unital completely positive maps on $M_n(\mathbb {C})$. 
\end{cor}
\smallskip 

{\it Proof}. Assume that $\Theta \in  Aut(M_n(\mathbb{C}))$ is given as an operational combination 
${\rm Ad}a \circ \Phi  + {\rm Ad} b \circ \Psi  = \Theta$ of  
$\Phi, \Psi  \in UCP(M_n(\mathbb{C}))$ with  an operational coefficient  $\{a, b\} \in FOP(M_n(\mathbb{C}))$.  
Then by Theorem 3.7, there exist  unitaries $v, w \in M_n(\mathbb{C})$ and a $\lambda \in (0,1)$ 
such that 
$a = \sqrt{\lambda} uv^*, \ b = \sqrt{1 - \lambda} uw^* \quad \text{and} \quad \Phi = {\rm Ad} v, \ \Psi(x) = {\rm Ad} w $ 
for a unitary $u$ with $\Theta = {\rm Ad} u$. 
These condition imply that $aa^* = \lambda 1_{M_n(\mathbb{C})}$, $bb^* = (1 - \lambda) 1_{M_n(\mathbb{C})}$ and that 
$\lambda^{-1} {\rm Ad} a \circ \Phi = (1 - \lambda)^{-1} {\rm Ad} b \circ \Psi =  \Phi$ 
so that $\Theta$ satisfies the condition of operational extreme points. 
\qed 

\end{document}